\input amstex 
\documentstyle{amsppt}
\input bull-ppt
\keyedby{bull413e/mhm}

\topmatter
\cvol{29}
\cvolyear{1993}
\cmonth{July}
\cyear{1993}
\cvolno{1}
\cpgs{1-13}
\title``Theoretical mathematics'':
Toward a
cultural \\synthesis of 
mathematics and theoretical physics
\endtitle
\author Arthur Jaffe and Frank Quinn \endauthor
\shorttitle{Theoretical Mathematics}
\address Harvard University, Cambridge, Massachusetts 
02138-2901\endaddress
\ml jaffe\@math.harvard.edu \endml
\address Virginia Polytechnical Institute and State 
University,
Blacksburg, Virginia 24061-0123\endaddress
\ml quinn\@math.vt.edu \endml
\date November 6, 1992\enddate
\thanks Both authors were partially supported by the 
National Science
Foundation\endthanks
\abstract Is speculative mathematics dangerous?  Recent 
interactions between
physics and mathematics pose the question with some force: 
 traditional
mathematical norms discourage speculation, but it is the 
fabric of theoretical
physics.  In practice there can be benefits, but there can 
also be unpleasant
and destructive consequences.  Serious caution is 
required, and the issue
should be considered before, rather than after, obvious 
damage  occurs.  With
the hazards carefully in mind, we propose a framework that 
should allow a
healthy and  positive role for speculation.\endabstract

\endtopmatter

\document


Modern mathematics is nearly characterized by the use of 
rigorous
proofs.  This practice, the result of literally thousands 
of years of
refinement, has brought to mathematics a clarity and 
reliability
unmatched by any other science. But it also makes 
mathematics slow and
difficult; it is arguably the most disciplined of human 
intellectual
activities.  

Groups and individuals within the mathematics community 
have from
time to time tried being less compulsive about details of 
arguments.  The
results have been mixed, and they have occasionally been 
disastrous.  
Yet today in certain areas there is again a trend toward 
basing
mathematics on intuitive reasoning without proof. To some 
extent this
is the old pattern of history being repeated by those 
unfamiliar with it.
But it also may be the beginning of fundamental changes in 
the way
mathematics is organized.  In either case, it is vital at 
this time to
reexamine the role of proofs in mathematical understanding 
and to
develop a constructive context for these trends. 

We begin with a discussion of physics, partially because 
some of the
current movement results from an interaction with 
theoretical physics 
and partly because it provides a useful model for 
potential sociological
realignments. We then turn to the history of mathematics 
for examples
illustrating the benefits and hazards of nonrigorous work. 
Finally, we 
suggest a framework which should allow different 
approaches to coexist.

\heading Theory and rigor\endheading

Typically, information about mathematical structures is 
achieved in two
stages.   First, intuitive insights are developed, 
conjectures are made, and
speculative outlines of justifications are suggested. Then 
the conjectures
and speculations are  corrected; they are made reliable by 
proving them. 
We use the term {\it theoretical\/} mathematics for the 
speculative and
intuitive  work; we refer to the proof-oriented phase as 
{\it rigorous\/}
mathematics.   We do not wish to get involved in a 
discussion of our
choice of terminology; this is not the central point of 
our article. 
However, as a point of departure, we briefly clarify our 
definitions of
theory and rigor. 

The initial stages of mathematical discovery---namely, the 
intuitive and
conjectural work, like theoretical work in the 
sciences---involves
speculations on the nature of reality beyond established 
knowledge.  Thus
we borrow our name ``theoretical'' from this use in 
physics.  There is an
older use of the word ``theoretical'' in mathematics, 
namely, to identify
``pure'' rather than applied mathematics; this is a usage 
from the past
which is no longer common and which we do not adopt.

Theoretical work requires correction, refinement, and 
validation through
experiment or proof.  Thus we claim that the role of 
rigorous proof in
mathematics is functionally analogous to the role of 
experiment in the
natural sciences.  This thesis may be unfamiliar but after 
reflection 
should be clear at least to mathematicians.  Proofs serve 
two main
purposes.  First, proofs provide a way to ensure the 
reliability of
mathematical claims, just as laboratory verification 
provides a check in
other sciences.  Second, the act of finding a proof often 
yields, as a
byproduct, new insights and unexpected new data, just as 
does work in
the laboratory.  

Mathematicians may have even better experimental access to
mathematical reality than the laboratory sciences have to 
physical reality.
This is the point of modeling: a physical phenomenon is 
approximated
by a mathematical model; then the model is studied 
precisely because it
is more accessible.   This accessibility also has had 
consequences
for mathematics on a social level.  Mathematics is much 
more finely
subdivided into subdisciplines than physics, because the 
methods have
permitted a deeper penetration into the subject matter. 

Although the use of proof in mathematics is functionally 
parallel to
experiment, we are not suggesting that proofs should be 
called
``experimental'' mathematics. There is already a 
well-established and
appropriate use of that term, namely, to refer to 
numerical calculations
and computer simulations as tests of mathematical 
concepts.  In fact,
results of computer experiments are frequently presented 
in a way we
would call theoretical: general conclusions are proposed 
on the basis of
test cases and comparisons. The conclusions are not 
completely reliable,
and an effort to provide real proofs would inevitably turn 
up exceptions
and limitations.  

\heading Division of labor\endheading

In physics  we have come to accept a division of labor 
between theorists
and experimentalists. But in fact only recently has the 
division become
so clear.  Until the beginning of the twentieth century 
there was
basically one community of physicists.  It was the ideal, 
and by and large
the practice, that the same people both speculated about 
theory and 
verified their speculations with laboratory experience.  

Certainly in Europe it had become clear by 1900 that a 
bifurcation had
occurred: there were sufficiently many physicists who 
concentrated solely
on the theoretical side of their work that one could 
identify two distinct
communities  [H].  This trend proceeded somewhat more 
slowly in the
United States. E.~C.~Kemble, who worked at Harvard in the 
area of
quantum theory, is generally regarded as the first 
American to obtain 
a doctorate in purely theoretical physics (though his 1917 
thesis
contained an experimental appendix).

Distinctness of the theoretical and experimental physics 
communities 
should not be confused with their independence. Theory is 
vital for
experimentalists to identify crucial tests and to 
interpret the data. 
Experiment is vital for theorists to correct and to guide 
their
speculations. Theoretical and experimental groups are 
unstable and
ineffective unless they occur in closely interacting pairs.

We contrast this division of labor in physics to the 
current situation in
mathematics.  It is still the ideal, and by and large the 
practice, that the
same people both speculate about mathematical structures 
and verify
their speculations through rigorous proofs.  In other 
words, the
mathematical community has not undergone a bifurcation  
into theoretical
and rigorous branches.  Even on the individual level we 
seldom
recognize theoretical mathematics as an appropriate 
principal activity.  

There have been attempts to divide mathematical efforts in 
this way, but
these attempts were for the most part unsuccessful.  Why?  
Is there
something about mathematics itself which invalidates the 
analogy with
physics and prevents such a bifurcation?  Is ``theoretical 
mathematics''
in the end an oxymoron?  Or were there flaws in past 
attempts which
doomed them but which might be avoided today?   

\heading New relations with physics\endheading

A new connection with physics is providing a good deal of 
the driving
force toward  speculation in mathematics.  Recently there 
has been a
flurry of mathemati\-cal-type activity in physics, under 
headings like
``string theory'', ``conformal field theory'', 
``topological quantum field
theory'', and ``quantum gravity''.  In large part this has 
been initiated by
individuals trained as theoretical high-energy physicists. 
 The most
celebrated and influential of these (though not the most 
problematic) is
Edward Witten. 

From a physical point of view much of this work has not 
yet matured
to the stage where observable predictions about nature 
have been made. 
Further, the work often concerns ``toy models'' designed 
to display only
analogs of real phenomena. And some of the parts which 
might apply to
the real world concern experimentally inaccessible events: 
 particles of
incredible energy, movement on the scale of the universe, 
or creation of
new universes.  

One result of the lack of predictions is that these 
physicists are cut off
from their presumptive experimental community; they have 
no source of
relevant physical facts to constrain and inspire their 
theorizing.   Since
progress comes from the interaction between theory and 
experiment,  a
theoretical group cannot exist long in isolation.  Indeed 
much of the
mainstream physics community regards these developments with
suspicion, because of their isolation from the so-called 
``real world''.  

But these physicists are not in fact isolated. They have 
found a new
``experimental community'':  mathematicians.  It is now 
mathematicians
who provide them with reliable new information about the 
structures they
study.  Often it is to mathematicians that they address 
their speculations
to stimulate new ``experimental'' work.  And the great 
successes are new
insights into mathematics, not into physics. What emerges 
is not a new
particle but a description of representations of the 
``monster'' sporadic
group using vertex operators in Kac-Moody algebras.  What 
is produced
is not a new physical field theory but a new view of 
polynomial
invariants of knots and links in 3-manifolds using Feynman 
path integrals
or representations of quantum groups.

These physicists are still working in the speculative and 
intuitive mode
of theoretical physics. Many have neither training for nor 
interest in 
rigor.  They are doing theoretical mathematics. Very senior
mathematicians have praised this work and have suggested 
it should be
emulated.  As a result, some mathematical followers are 
moving toward
a more speculative mode.

One could conclude from this description that parts of 
mathematics have
already been propelled through the bifurcation.  These 
areas have
suddenly acquired a fully functioning theoretical 
community, the
theoretical physicists, and traditional mathematicians 
have become the
partner community of rigorous verifiers.

However, this has happened without the evolution of the 
community
norms and standards for behavior which are required to 
make the new
structure stable.  Without rapid development and adoption 
of such
``family values'' the new relationship between mathematics 
and physics
may well collapse.  Physicists will go back to their 
traditional partners;
rigorous mathematicians will be left with a mess to clean 
up; and
mathematicians lured into a more theoretical mode by the 
physicists'
example will be ignored as a result of the backlash.  To 
understand what
is involved, we start with a look at the past.   

\heading Old relations with physics\endheading

The rich interplay between mathematics and physics 
predates even their
recognition as separate subjects. The mathematical work 
that in some
sense straddles the boundaries between the two is commonly 
referred to
as {\it mathematical physics\/}, though a precise 
definition is probably
impossible.  In particular over the past ninety years or 
so, a school of
mathematical physicists emerged around such persons as  
D.~Hilbert, 
F.~Klein, H.~Poincar\'e, M.~Born, and later H.~Weyl, 
J.~von~Neumann,
E.~P.~Wigner, M.~Kac, A.~S.~Wightman, R.~Jost, and 
R.~Haag---persons 
with training both in physics and in mathematics.  These 
people
often worked on questions motivated by physics, but they 
retained the
traditions and the values of mathematics. 

Results emerging from this school have at various times 
been relevant
both for mathematics and for physics.  Their development 
has proceeded
at a deliberate pace with the accumulation of conclusions 
of long-term
interest.  A few (out of many) recent examples include 
work on the
existence of quantum field theory and on its compatibility 
with relativity;
the work of Lieb and of Baxter on lattice models and 
related transfer
matrices;  the work of Schoen and Yau on the positive 
energy theorem
in relativity and its relation to the geometry of minimal 
surfaces; the
work of Ruelle and others on dynamical systems and 
turbulence;  the
operator algebra approach to local quantum theory; and 
Connes's early
interest in physics, from which emerged his mathematical 
work on
factors and later the foundation of noncommutative 
geometry.   

The main point here is what has {\it not} happened. Work 
from this
school is characterized by standards of scholarship and by 
knowledge of
the literature consistent with the best traditions of 
mathematics.  There
is no ambiguity about definitions, the formulation of 
claims, or proofs of
theorems. It is traditional rigorous mathematics. Even in 
this close
proximity to physics the traditional values of mathematics 
have been
retained. Mathematical physicists do have access to a 
vast, rich body of
speculation 
by theoretical physicists. But these speculations have
traditionally been addressed to physicists, not 
mathematicians.  

Theoretical physics and mathematical physics have rather 
different
cultures, and there is often a tension between them. 
Theoretical work in
physics does not need to contain verification or proof, as 
contact with
reality can be left to experiment.  Thus the sociology of 
physics tends to
denigrate proof as an unnecessary part of the theoretical 
process. 
Richard Feynman used to delight in teasing mathematicians 
about their
reluctance to use methods  that ``worked'' but  that could 
not be
rigorously justified \cite{F, G2}.   He felt it was quite 
satisfactory
to test mathematical statements by verifying a few 
well-chosen cases.   

On the mathematical side, E.~J.~McShane once likened the 
reasoning in
a ``physical argument'' to that of ``the woman who could 
trace her
ancestry to William the Conqueror with only two gaps'',  
and this was
typical of the mathematical attitude. There was an 
understanding that the
development must proceed at an internally appropriate pace 
and not be
stampeded by the extent of the physicists' vision.   
McShane and most
of the mathematical physics community rejected the 
free-wheeling,
theoretical approach as inappropriate. 

A relevant observation is that most theoretical physicists 
are quite
respectful of their experimental counterparts.  Relations 
between physics
and mathematics would be considerably easier if physicists 
would
recognize mathematicians as ``intellectual 
experimentalists'' rather than
think of them disdainfully as  uselessly compulsive 
theorists.  The typical
attitude of physicists toward mathematics is illustrated 
by a passage from
a book of P.~W.~Anderson, ``We are talking here about 
theoretical
physics, and therefore of course mathematical rigor is 
irrelevant and
impossible.'' \footnote{
This is attributed by Anderson to
Landau \cite{A, p.\ 132}. Anderson continues, ``This is 
not quite so,
but it is  very close to it.'' However, he revealingly 
remembered the
passage  incorrectly; it reads  \cite{LL},  ``No attempt 
has been made at
mathematical rigor in the treatment, since this is anyhow 
illusory in
theoretical physics,$\ldots$.''} In fact it is exactly as 
relevant and
possible as experimental data, and like data should  
be used whenever available.  Nevertheless, students in 
physics are
generally indoctrinated with antimathematical notions; and 
if they
become involved in mathematical questions, they tend not 
only to be
theoretical but often to deny that their work is incomplete.

Not all of mathematical physics has been as clearly 
mathematical as that
described above.  For example, the work of the German 
mathematical
physicist K.~Symanzik was mostly theoretical. However, he 
was very
careful not to make unwarranted claims of mathematical 
rigor.  In fact,
he made a serious attempt in 1968 to establish an 
important part of his
theoretical program on a rigorous level in collaboration 
with the
mathematician S. R. S.~Varadhan.  Some years afterward, 
this was
achieved by E. Nelson, Osterwalder, Schrader, and others 
and had
far-reaching consequences relating quantum theory, 
probability, and
statistical mechanics.  

We also mention two people who have been addressing largely 
mathematical questions from an almost entirely theoretical 
point of view
for some time---namely, B.~Mandelbrot and M.~Feigenbaum.  
See
\cite{G} for a popular account and \cite {Kr}
for an  expression of the mathematical discomfort with 
this activity. 

For the most part the mainstream of mathematical physics 
has rejected
purely theoretical work as a valid mathematical style. We 
observe,
however, that the mathematicians involved in the ``new 
relations'' with
physics are different from the traditional mathematical 
physicists. 
Geometers, topologists, and persons in representation 
theory have begun
to talk with  physicists. These mathematicians are unused 
to dealing with
the difference in cultures and for the most part do not 
recognize parallels
in their own collective experience that would sensitize 
them to the
hazards of theoretical work.  This suggests a question: As 
they gain more
experience, will these mathematicians also reject pure 
theorizing? 

\heading Success stories\endheading

We turn to mathematical encounters with theoretical work 
and begin
with the positive side. The posing of conjectures is the 
most obvious
mathematical activity that does not involve proof.  
Conjectures range
from brilliant to boring, from impossible to obvious.  
They are filtered
by the interest they inspire rather than by editors and 
referees. The better
ones have inspired the development of whole fields. Some 
of the most
famous examples are the Riemann Hypothesis, Fermat's 
``last theorem'',
and the Poincar\'e conjecture.  The Hilbert problem list, 
of amazing
breadth and depth, has been very influential in the 
development of
mathematics in this century.  Other examples are the Adams 
conjecture
in topology, the several Questions of Serre, the Novikov 
conjecture, and
the Wightman axioms for quantum field theory.

Some conjectures are accompanied by technical details or a 
proposal for
a proof.  For example, the ``Weil conjectures'' outlined 
an approach to
a $p$-adic analog of the Riemann  hypothesis.  The 
implementation of
this program by Grothendieck and Deligne was celebrated as 
a major
achievement for modern algebraic geometry.  Similarly 
celebrated was
Falting's proof of the Mordell conjecture, part of a 
program to approach
the Fermat ``theorem''. The ``Langlands program'' for 
understanding
automorphic forms has been a major stimulus to that field, 
and ``Mori's
program'' for the investigation of algebraic three-folds 
has invigorated that
area.   The classification of finite simple groups was 
achieved by
implementing a program developed by Gorenstein. 

These examples share the characteristic that they were 
explicitly
speculative when formulated (or at least quickly 
recognized as
speculative, as for example with Fermat).  They 
represented a goal to
work toward, and primary credit for the achievement was 
clearly to be
assigned to the person who found a proof.   

Another type of mathematical work is intermediate between 
traditional
and theoretical.  It proceeds in the way, ``If A is true, 
then X, Y, and Z
follow'', or ``If A, then it is reasonable to conjecture 
R, S, and T.''  In
this case ``A'' may be a major outstanding mathematical 
conjecture,  the
Riemann Hypothesis for instance.  A striking recent 
example of work in
this style occurs in the theory of motifs, and recently an 
entire
S\'eminaire Bourbaki was devoted to the evolution of work 
in number
theory and algebraic geometry based on conjectures of 
Deligne and of
Beilinson \cite{Fo}.  It is interesting that the Bourbaki, 
once the bastion
of the most conservative traditional mathematics, now 
entertains
pyramids of conjectures.  
 
The importance of large-scale, goal-formulating work 
(necessarily
theoretical) is growing.  We are in an age of big science, 
and
mathematics is not an exception. The classification of 
finite simple
groups, for instance, is estimated to occupy 15,000 
journal pages! Other
sciences have responded to this trend by forming large 
formal
collaborative groups.   

The National Science Foundation and government agencies in
other countries have tried to nudge mathematics in the 
direction of
collaborative and interdisciplinary work. But one must 
recognize that
large projects in mathematics are not centered around a 
grant, a
technique, or a machine.  They are undertaken by informal 
communities
nucleated by a visionary theoretical program. In the 
finite simple group
effort, the program was developed and coordinated by 
D.~Gorenstein,
who parcelled out pieces of the puzzle to an informal 
``team'' working
on the problem. Future growth of such large-scale 
mathematical activity
can only occur with the evolution of more such visionary 
programs.   

\heading Cautionary tales\endheading

Most of the experiences with theoretical mathematics have  
been less
positive than those described above. This has been 
particularly true
when incorrect or speculative material is presented as 
known and reliable,
and credit is claimed by the perpetrator.  Sometimes this 
is an ``honest
mistake'', sometimes the result of nonstandard conceptions 
of what
constitutes proof.  Straightforward mistakes are less 
harmful.  For
example, the fundamental ``Dehn lemma'' on two-disks in 
three-manifolds was
presented in 1910.  An error was found, and by the time it 
was proved
(by Papakyriakopolos in 1957) it was recognized as an 
important
conjecture. 

Weak standards of proof cause more difficulty.  In the 
eighteenth
century, casual reasoning led to a plague of problems in 
analysis
concerning issues like convergence of series and uniform 
convergence of
functions.  Rigor was introduced as the antidote.  It was 
adopted over
the objections of some theorists in time to avoid major 
damage.

More recently in this century the ``Italian school'' of 
algebraic
geometry did not avoid major damage:  it collapsed after a 
generation of
brilliant speculation.  See \cite{EH, K} for discussions 
of the
difficulties and the long recovery.  In 1946 the subject 
was still regarded
with such suspicion that Weil felt he had to defend his 
interest in it;
see the introduction to \cite{W}.

Algebraic and differential topology have had several 
episodes of
excessively theoretical work.  In his history \cite{D}, 
Dieudonn\'e dates
the beginning of the field to Poincar\'e's {\it Analysis 
Situs} in 1895. 
This ``fascinating and exasperating paper'' was extremely 
intuitive.  In
spite of its obvious importance it took fifteen or twenty 
years for real
development to begin. Dieudonn\'e expresses surprise at 
this slow start
\cite{D, p.\ 36}, but it seems an almost inevitable 
corollary of how it
began: Poincar\'e claimed too much, proved too little, and 
his ``reckless''
methods could not be imitated.  The result was a dead area 
which had to
be sorted out before it could take off. 

Dieudonn\'e suggests that casual reasoning is a childhood 
disease of
mathematical areas and says, ``\<$\dotsc$\<after 
1910\<$\dotsc$\<uniform 
standards of
what constitutes a correct proof became universally 
accepted in
topology\<$\dotsc$\<this standard has remained {\it 
unchanged} ever since.''
But in fact there have been many further episodes.  Ren\'e 
Thom's early
work on differentiable manifolds, for which he received 
the Fields medal,
was brilliant and generally solid. Later work on 
singularities was not so
firm.  His claim of $C^{\infty}$ density of topologically 
stable maps was
supported by a detailed but incomplete outline,  which was 
later repaired
by John Mather.  Thom went on to propose the use of 
``catastrophe
theory'', founded on singularities, to explain forms of 
physical
phenomena.  The application was mathematically 
theoretical, and its
popularization, particularly by E.~C.~Zeeman, turned out 
to be physically
controversial. 

The early work of Dennis Sullivan provides another 
example. After a
solid beginning, in the 1970s he launched into a brilliant 
and highly
acclaimed but ``theoretical'' exploration of the topology 
of manifolds. 
Details were weak, and serious efforts to fill them in got 
bogged down. 
Sullivan himself changed fields and returned to a more 
rigorous
approach.  This field still seems to have more than its 
share of fuzzy
proofs. 

William Thurston's ``geometrization theorem'' concerning 
structures
on 
Haken three-manifolds is another often-cited example. A 
grand insight
delivered with beautiful but insufficient hints, the proof 
was never fully
published.  For many investigators this unredeemed claim 
became a
roadblock rather than an inspiration.

In these examples, as with Poincar\'e, the insights 
proposed seem to be 
on target.   There are certainly cases in which the 
theoretical insights
were also flawed.   The point is that even in the best 
cases there were
unpleasant side effects that might have been avoided. We 
also see that
Witten, in giving a heuristic description of an extension 
of the Jones
polynomial \cite{Wi}, was continuing in a long and 
problematic tradition
even within topology. 

Some areas in the Russian school of mathematics have 
extensive
traditions of theoretical work, usually conducted through 
premature
research announcements. From the numerous possible 
examples we
mention only two.   The first is concerned with  the 
perturbation theory
of integrable Hamiltonian systems with phase space 
foliated by invariant
tori.  In 1954 Kolmogorov announced that tori with 
nonresonant
frequencies survive a perturbation and gave an outline of 
an argument.
In retrospect it may be seen that this outline does touch 
on the major
ideas necessary, but it was generally considered  
insufficient to allow
reconstruction of a proof. Complete proofs were achieved 
by Arnold in
1959 for the analytic case and by Moser in 1962 for the 
smooth case. 

The second example is one in which one of us became 
personally
involved while working on trying to establish a widely 
conjectured result
that phase transitions occur in (relativistic) quantum 
field theory.  
In 1973 the respected mathematicians Dobrushin and Minlos 
published
an announcement of that result.  Two years later when no 
indication
had come from the Russians of a proof, Glimm, Jaffe, and 
Spencer
resumed their work on the problem and eventually gave two 
different
proofs.  A couple of years after that Dobrushin and Minlos 
published a
retraction of their original announcement.

\heading The problems\endheading

There are patterns in the problems encountered in these 
examples. We
list some and then discuss them in more detail.  

(1) Theoretical work, if taken too far, goes astray 
because it lacks
the feedback and corrections provided by rigorous proof.

(2) Further work is discouraged and confused by uncertainty
about which parts are reliable. 

(3) A dead area is often created when full credit is 
claimed by
vigorous theorizers:  there is little incentive for 
cleaning up the debris
that blocks further progress.  

(4) Students and young researchers are misled. 



\noindent The first problem often overtakes would-be 
mathematical
theorizers, particularly when they are unwilling to 
acknowledge that their
work is uncertain and incomplete.  Even in theoretical 
physics where
there is an awareness of this possibility, knowing when to 
stop is a subtle
and difficult skill. Errant theorizing damages the 
credibility of the
theorist and may also damage the field through the 
mechanism identified
in the second problem.

The second problem has to do with uncertainty of the 
literature. In
comparison with other sciences, the primary mathematical 
literature is
extraordinarily reliable.  Papers in refereed core 
mathematics journals 
are nearly always sound, and this permits steady and 
efficient advance. 
Only a small pollution of serious errors would force 
mathematicians to
invest a great deal more time and energy in checking 
published material
than they do now.  The advantages of a reliable literature 
are so profound
that we suspect this is the primary force driving 
mathematics toward
rigor.

When reliability of a literature is uncertain, the issue 
must be addressed. 
Often ``rules of thumb'' are used.  
For example, mathematicians presume that papers in physics 
journals are
theoretical.  This extends to a suspicion of 
mathema\-tical-physics journals,
where the papers are generally reliable (though with 
dangerous
exceptions).  Another widely applied criterion is that 
anything using
functional integrals must be speculative.  One of us has 
remarked on the
difficulties this causes mathematicians trying to use 
solid instances of the
technique \cite{J}.  These kinds of rules are 
unsatisfactory, as is the {\it
caveat emptor} approach of letting each paper be judged 
for itself. 
Proponents of this latter view cite Witten's papers as 
successful
examples.  But a few  instances can be handled; it is 
large numbers that
are a disaster.  Also, it is a common rule of thumb now to 
regard any
paper by Witten as theoretical. This short-changes 
Witten's work but
illustrates the ``better-safe-than-sorry'' approach 
mainstream mathematics
tends to take when questions arise.  

This unreliability is certainly a problem in theoretical 
physics, where the
primary literature often becomes so irrelevant that it is 
abandoned
wholesale.  I.~ M.~Singer has compared the physics 
literature to a
blackboard that must be periodically erased. Physicists 
traditionally
obtain  much less benefit from the historical background 
of a problem,
and they are less apt to search the literature. The 
citation half-life of
physics papers is much shorter than in mathematics.   

The ``dead area'' problem concerns credit and rewards. 
Mathematical
researchers traditionally do not give credit twice for the 
same results. But
this means that when a theorist claims credit, it is 
difficult for rigorous
workers to justify the investment of labor required to 
make it reliable. 
There is a big difference between ``filling in the details 
of a theorem by
$X$'' and ``verifying a conjecture of $X$''.   Rigorous 
mathematicians
tend to flee the shadow of a big claim. The pattern is 
that the missing
work is filled in, often much later, using techniques and 
corollaries of
work on separate topics for which uncontested credit is 
available. 

Finally, on the last point  most successful theorizers (at 
least in
mathematics) have a solid background in disciplined work, 
which is the
source of their intuition and taste.  Most students who 
try to dive directly
into the heady world of theory without such a background are
unsuccessful.  Failure to distinguish between the two 
types of activity can
lead students to try to emulate the more glamorous and 
less disciplined
aspects and to end up unable to do more than manipulate 
jargon.

Mathematicians tend to focus on intellectual content and 
neglect the
importance of social issues and the community.  But we 
{\it are} a
community and often form opinions even on technical issues 
by social
interactions rather than directly from the literature.  
Socially accepted
conventions are vital in our understanding of what we 
read.  Behavior is
important, and the community of mathematicians is 
vulnerable to damage
from inappropriate behavior.     

\heading Prescriptions\endheading

The mathematical community has evolved strict standards of 
proof and 
norms  that discourage speculation.  These are protective 
mechanisms
that guard against the more destructive consequences of 
speculation; they
embody the collective mathematical experience that the 
disadvantages
outweigh the advantages.  On the other hand, we have seen 
that
speculation, if properly undertaken, can be profoundly 
beneficial. 
Perhaps a more conscious and controlled approach that
would allow us to reap the benefits but avoid the dangers 
is possible. The need to 
find a constructive response to the new influences from 
theoretical
physics presents us with both an important test case and 
an opportunity. 

Mathematicians should be more receptive to theoretical 
material but with
safeguards and a strict honesty.  The safeguards we 
propose are not new;
they are essentially the traditional practices associated 
with conjectures.
However, a better appreciation of their function and 
significance is
necessary, and they should be applied more widely and more 
uniformly. 
Collectively, our proposals could be regarded as measures 
to ensure
``truth in advertising''. 
$$\topctrblock
Theoretical work should be explicitly
acknowledged as theoretical and incomplete; in particular, 
a major share
of credit for the final result must be reserved for the 
rigorous work that
validates it.
\endblock$$
\noindent This can make the difference between a dead area 
and a living
industry. Theoreticians should recognize that in the long 
run the success
of their work is dependent on the work of a companion 
rigorous
community; they should honor and nurture it when possible. 
Certainly in
physics the community assigns basic credit for discovery 
to successful
experimental investigations; they are not regarded merely 
as verifying
small details in the web of theoretical insight.  On the 
individual level
mathematical authors should make a choice:  either they 
provide
complete proofs, or they should agree that their work is 
incomplete and
that essential credit will be shared.  Referees and 
editors should enforce
this distinction, and it should be included in the 
education of students.

The other suggestions are concerned with the integrity of 
the
mathematical literature.  It has always been acceptable to 
state a
conjecture in a paper, and occasionally papers have been 
published  that
are entirely theoretical.  The key issue is to identify 
the theoretical
material clearly.  

$$\topctrblock Within a paper, standard
nomenclature should prevail:  in theoretical material, a 
word like
``conjecture'' should replace ``theorem''; a word like 
``predict'' should
replace ``show'' or ``construct''; and expressions such as 
``motivation''
or ``supporting argument'' should replace ``proof''. 
Ideally the title and
abstract should contain a word like ``theoretical'', 
``speculative'', or
``conjectural''. 
\endblock$$

\noindent The objective is to have flags indicating the 
nature of the work
to readers.  A flag in the title of a theoretical paper  
would appear in a
citation, which would help limit second-hand problems. 
Theoretical work
should be cited as a source of inspiration or to justify 
significance or in
supporting arguments in other theoretical developments.  
Citing a
theoretical paper for a structural ingredient of a 
supposedly rigorous
proof must be handled with care, and a flag in the title 
would indicate
when such care is needed. 

Research announcements pose particular problems. Some 
announcements
are simply summaries of work  that the author has 
completed and written
down. In such announcements, the language of theorem and 
proof is
appropriate. Others describe outcomes of arguments which 
have not been
worked out in detail, and sometimes they contain leaps of 
faith which
require years of effort to bridge. In these cases the 
traditional language
misrepresents the work and is not appropriate;  such 
announcements really
should be identified as theoretical.  The analysis above 
implies that
because of this, publication of announcements may be 
unwholesome or
even damaging to the fabric of mathematics.
Announcements do have valid functions, for instance, 
establishing a claim
to priority and alerting others to new results and useful 
techniques.  We
therefore seek guidelines, as with theoretical work in 
general, which will
permit beneficial uses but limit the potential damage. The 
key issue
seems to be incorporation into the literature.  
One solution is:
$$\topctrblock Research announcements  should not be
published, except as summaries of full versions  that have 
been accepted
for publication.   Citations of unpublished work should 
clearly distinguish
between announcements and complete preprints. 
\endblock$$

\noindent 
In this age of copying machines and electronic bulletin 
boards it is
possible to distribute information widely without formally 
publishing it.
``Cross-cultural'' preliminary results of wide interest 
could be described
informally in ``news column'' format in publications like 
the
{\it Mathematical Intelligencer\/} or the 
{\it Notices of the American Mathematical
Society\/}. The lack of formal publication, therefore, 
should be at most a
minor disadvantage. We observe that this discussion also 
underscores the
importance of maintaining the distinction between formal 
(refereed)
publication and a posting on a bulletin board. Maintaining 
this distinction
may be one of the greatest challenges facing the 
development of serious
electronic journals.

If these safeguards are carefully followed, then it would 
be reasonable for
any mathematics journal to consider theoretical articles 
for publication. 
Stimulating articles with incomplete proofs might be 
offered publication
(after appropriate word changes) as theoretical papers.  
Theoretical
mathematics journals  might be appropriate.  But care is 
required.
Without honesty and caution by authors, editors, and 
referees this would 
simply lead to the reintroduction of problems  that have 
been painfully
and repeatedly purged over the years.

Our analysis suggests that the bifurcation of mathematics 
into theoretical 
and rigorous communities has partially begun but has been 
inhibited by
the consequences of improper speculation.  Will it 
continue?  
Probably it will in any case, but it should evolve faster 
and less painfully
if the safeguards are adopted.  In the classical areas it 
will be slower
because intelligent speculation must be based on a mastery 
of technical
detail in previous proofs.  In these areas the framework 
for speculation
is more likely to provide a constructive outlet for 
nominally rigorous
individuals, whose inspirations exceed their capacity for 
rigorous proof. 
Other areas, particularly those involving computer 
simulations, are
different in that the generation and analysis of data is 
quite a distinct
activity from the construction of proofs. In some of these 
areas
specialized theorists can already be identified and are 
likely to
proliferate.

In any case, the proposed framework provides a context for 
interactions
between mathematicians and theoretical physicists.  
Whether or not they
become a permanent fixture in the mathematical community, 
physicists
can be welcomed as ``theoretical mathematicians'' rather 
than rejected
as incompetent traditional mathematicians.  

\heading Summary\endheading

At times speculations have energized development in 
mathematics; at
other times they have inhibited it.  This is because  
theory and proof are
not just ``different'' in a neutral way.  In particular, 
the failure to
distinguish carefully between the two can cause damage  
both to the
community of mathematics and to the mathematics 
literature.  One might
say that it is mathematically unethical not to maintain 
the distinctions
between casual reasoning and proof.  However, we have 
described
practices and guidelines which, if carefully implemented, 
should give a
positive context for speculation in mathematics. 

\heading Acknowledgment \endheading

We wish to thank many colleagues who have made helpful 
suggestions about this
paper. The first author thanks the John S. Guggenheim 
Foundation for a
fellowship.   Both authors were partially supported by the 
National Science
Foundation.


\Refs\widestnumber\key{EH}

\ref\key A\by P.~W.~Anderson\book Concepts in solids
\publaddr New York
\publ
W.~A.~Benjamin, Inc\yr1964\endref

\ref\key D\by J.  Dieudonn\'e\book A history of algebraic 
and
differential topology {\rm 1900--1960}
\publaddr Basel 
\publ Birkh\"auser\yr1988\endref

\ref\key EH\by D.  Eisenbud and J. Harris \paper Progress 
in the theory
of complex algebraic curves\jour Bull. Amer. Math. Soc. 
\vol21\yr1989\page205\endref

\ref\key F\by R. P. Feynman\book Surely you\RM're joking 
Mr. Feynman\/\RM: 
adventures of a curious character\publ W. W. Norton
\publaddr New York\yr 1985\endref 

\ref\key Fo\by J.-M. Fontaine\paper Valeurs sp\'eciales 
des fonctions $L$
des motifs \inbook S\'eminaire Bourbaki, Expos\'e 751, 
F\'evrier
1992\pages1--45\endref 

\ref\key G1 \book Chaos\,\RM: making a new science\by J. 
Gleick\publ Viking
Penguin Inc.\publaddr New York\yr 1987\endref 

\ref\key G2 \book Genius\,\RM: the life and science of 
Richard Feynman \bysame
\publ Pantheon \publaddr New York\yr 1992\endref 

\ref\key H\by G. Holton\publ private communication\endref 

\ref\key J\by A.  Jaffe \paper Mathematics motivated by 
physics\inbook
Proc. Sympos. Pure Math.\publ Amer. Math. Soc.
\publaddr Providence, RI
\vol50 \yr1990\pages137--150\endref 

\ref\key K\by J.  Kollar \paper The structure of algebraic 
threefolds\,\/\RM:  an
introduction to Mori\RM's program\jour Bull. Amer. Math. 
Soc.
\vol17\yr1987\page 211\endref 

\ref\key Kr\by S. G. Krantz\paper Fractal geometry\jour 
Math. 
Intelligencer \vol 11\pages 12--16\yr 1989\endref

\ref\key LL\book Statistical physics \by L. Landau and E. 
Lifshitz \publ
Oxford Univ. Press
\publaddr London \yr 1938 \endref 

\ref\key M\by R. McCormmach, ed. \book Historical studies 
in the
physical sciences\publ Princeton Univ. Press
\publaddr Princeton, NJ \yr1975 \endref

\ref\key W\by A. Weil\book Foundations of algebraic 
geometry\publ
Amer. Math. Soc. \publaddr Providence, RI \yr 1946\endref

\ref\key Wi\by  E. H. Witten \paper Quantum field theory 
and the  Jones
polynomial\jour Comm. Math. Phys.\vol 121\yr 1989 \pages
351--399\endref

\endRefs
\enddocument